\documentclass[12pt]{article}
\usepackage{mathrsfs}
\usepackage{epic,eepic,epsf,epsfig}
\usepackage{amsfonts,srcltx,mathrsfs}
\usepackage{multirow}
\usepackage{amsmath}
\usepackage{amssymb}
\usepackage{amsbsy}
\usepackage{graphicx}
\usepackage{amsfonts}
\usepackage{color}
\usepackage{setspace}

\newtheorem{lem}{Lemma}[section]%
\newtheorem{theorem}[lem]{Theorem}%

\def\nd{\mathrel{\bigm|\kern-.7em/}}

\def\f{\noindent}

\def\P\GammaL{\hbox{\rm P\GammaL}}

\def\mod{\hbox{\rm mod }}

\begin{document}
\title{Some lemmas on spectral radius of graphs: including an application}

\footnotetext{E-mails: zhangwq@pku.edu.cn}

\author{Wenqian Zhang\\
{\small School of Mathematics and Statistics, Shandong University of Technology}\\
{\small Zibo, Shandong 255000, P.R. China}}
\date{}
\maketitle

\begin{abstract}
For a graph $G$, the spectral radius $\rho(G)$ of $G$ is the largest eigenvalue of its adjacency matrix. In this paper, we give three lammas on $\rho(G)$ when $G$ contains a spanning complete bipartite graph. Using these lemmas and typical spectral method, we characterized the unique extremal graph with the maximum spectral radius among all planar graphs of large order $n$ without a cycle of length $\ell$, where $5\leq \ell\leq n$.

\bigskip

\f {\bf Keywords:} spectral radius; planar graph; cycle; walk.

\end{abstract}

\baselineskip 17 pt

\section{Introduction}

All graphs considered in this paper are finite, undirected and simple.
For a graph $G$, let $\overline{G}$ be its complement. The vertex set and edge set of $G$ are denoted by $V(G)$ and $E(G)$, respectively. For two vertices $u,v$, let $d_{G}(u,v)$ denote the distance between them (i.e., the minimum length of a path connecting them). We say $u\sim v$, if $u$ and $v$ are adjacent in $G$.  For a certain integer $n$, let $K_{n},P_{n}$ and $C_{n}$ be the complete graph, the path and the cycle of order $n$, respectively.  For any $\ell\geq2$ graphs $G_{1},G_{2},...,G_{\ell}$, let $\cup_{1\leq i\leq \ell}G_{i}$ be the disjoint union of them. In particular, let $tG=\cup_{1\leq i\leq t}G$ for any integer $t\geq1$. Let $G_{1}\vee G_{2}$ be the join of $G_{1}$ and $G_{2}$, obtained from $G_{1}\cup G_{2}$ by connecting all vertices of $G_{1}$ to all vertices of $G_{2}$. For any terminology used but not defined here, one may refer to \cite{CRS}.

Let $G$ be a graph with vertices $v_{1},v_{2},...,v_{n}$. The {\em adjacency matrix} $A(G)$ of $G$ is an $n\times n$ matrix $(a_{ij})$, where $a_{ij}=1$ if  $v_{i}\sim v_{j}$, and $a_{ij}=0$ otherwise. The {\em spectral radius} $\rho(G)$ of $G$ is the largest eigenvalue of its adjacency matrix. By the well-known Perron--Frobenius theorem,  $\rho(G)$ has a non-negative eigenvector.  A non-negative eigenvector corresponding to $\rho(G)$ is called a {\em Perron vector} of $G$. If $G$ is connected, any Perron vector of $G$ has positive entries. The spectral radius of a graph is  closely related to its convergence property. In this paper, we first give the following three lammas on spectral radius of graphs. (Let $P_{0}$ denote the empty graph.)

 \begin{lem}\label{1-path}
For integers $n_{1}\geq n_{2}+2\geq3$, let $$G_{1}=H\vee(P_{n_{1}}\cup P_{n_{2}}\cup T)$$ and $$G_{2}=H\vee(P_{n_{1}-1}\cup P_{n_{2}+1}\cup T),$$ where $H,T$ are two graphs with $|H|\geq1$ and $|T|\geq0$. Then $\rho(G_{1})>\rho(G_{2})$.
\end{lem}

 \begin{lem}\label{2-path}
For integers $n_{1}\geq n_{2}\geq n_{3}\geq n_{4}\geq1$ and $n_{1}\geq n_{2}+2$, let $$G_{1}=H\vee(P_{n_{1}}\cup P_{n_{2}}\cup P_{n_{3}}\cup P_{n_{4}}\cup T)$$ and $$G_{2}=H\vee(P_{n_{1}-1}\cup P_{n_{2}+1}\cup P_{n_{3}+1}\cup P_{n_{4}-1}\cup T),$$ where $H,T$ are two graphs with $|H|\geq1$ and $|T|+\sum_{1\leq s\leq4}n_{s}\geq2000$. Then $\rho(G_{1})<\rho(G_{2})$.
\end{lem}

 \begin{lem}\label{3-path}
For integers $n_{1}\geq n_{2}\geq n_{3}\geq n_{4}\geq n_{5}\geq1$ and $n_{1}\geq n_{2}+3$, let $$G_{1}=H\vee(P_{n_{1}}\cup P_{n_{2}}\cup P_{n_{3}}\cup P_{n_{4}}\cup P_{n_{5}}\cup T)$$ and $$G_{2}=H\vee(P_{n_{1}-2}\cup P_{n_{2}+1}\cup P_{n_{3}+1}\cup P_{n_{4}+1}\cup P_{n_{5}-1}\cup T),$$ where $H,T$ are two graphs with $|H|\geq1$ and $|T|+\sum_{1\leq s\leq5}n_{s}\geq6500$. Then $\rho(G_{1})<\rho(G_{2})$.
\end{lem}

 A graph is called planar, if it can be drawn in the plane without crossed edges.  Tait and Tobin \cite{TT} proved that $K_{2}\vee P_{n-2}$ is the only extremal graph with the maximum spectral radius among all planar graphs of large order $n$. This (asymptotically) solves a conjecture proposed by Boots and Royle \cite{BR}, and
independently by Cao and Vince \cite{CV}.  Recently, the researchers are also interested in determining the extremal graphs with the maximum spectral radius among all $C_{\ell}$-free planar graphs of large order $n$, where $3\leq \ell\leq n$. The case $\ell=4$ can be settled from the results of  Nikiforov \cite{N}, and  Zhai and Wang \cite{ZW}. Very recently, Fang, Lin and Shi \cite{FLS} solved the case when $\ell=3$ and $5\leq \ell\ll n$. Xu, Lin and Fang \cite{XLF} solved the case when $\frac{\ell}{n}$ is close to 1. They \cite{XLF} proposed the general case as a problem. Armed with Lemma \ref{1-path} and Lemma \ref{2-path}, using typical spectral method we can solve this problem. That is, we have the following Theorem \ref{cycle-planar}.

 For a finite family of graphs $\mathcal{G}$, let ${\rm SPEX}(\mathcal{G})$ be the set of all extremal graphs  with the maximum spectral radius in $\mathcal{G}$. 
 
 \begin{theorem}\label{cycle-planar}
For integers $n\geq1.8\times10^{17}$ and $5\leq\ell\leq n$, let $\mathcal{G}_{n,\ell}$ be the set of all $C_{\ell}$-free planar graphs of order $n$. Then we have the following conclusions.\\
$(i)$ For $5\leq \ell<\frac{2n+5}{3}$, ${\rm SPEX}(\mathcal{G}_{n,\ell})=
\left\{K_{2}\vee(P_{\lceil\frac{\ell-3}{2}\rceil}\cup (a+1)P_{\lfloor\frac{\ell-3}{2}\rfloor}\cup P_{b})\right\}$, where $n-\ell+1=a\lfloor\frac{\ell-3}{2}\rfloor+b$ with $0\leq b<\lfloor\frac{\ell-3}{2}\rfloor$.\\
$(ii)$ For $\frac{2n+5}{3}\leq\ell\leq n$, ${\rm SPEX}(\mathcal{G}_{n,\ell})=
\left\{K_{2}\vee(P_{2\ell-n-4}\cup2P_{n-\ell+1})\right\}$.
 \end{theorem}

 The rest of the paper is organized as follows. In Section 2, we include several lemmas, which will be used later. The proofs of Lemmas \ref{1-path}, \ref{2-path} and \ref{3-path} are given in Section 3. In Section 4, we give the proof of Theorem \ref{cycle-planar}.

\medskip

 \f{\bf Note.} In this third version, we correct some calculation errors in Lemma \ref{w-evaluation}.

\section{Preliminaries}

Let $G$ be a graph. For an integer $\ell\geq1$, a {\em walk} of length $\ell$ in $G$ is an ordered sequence of vertices $v_{0}v_{1}\cdots v_{\ell}$, such that $v_{i}\sim v_{i+1}$ for any $0\leq i\leq \ell-1$. The vertex $v_{0}$ is called the starting vertex of the walk.  For any $u\in V(G)$ and $\ell\geq1$, let $w^{\ell}_{G}(u)$ be the number of walks of length $\ell$ starting at $u$ in $G$. Let $W^{\ell}(G)=\sum_{v\in V(G)}w^{\ell}_{G}(u)$.
The following spectral formula  was given in \cite{Z}.

\begin{theorem}{\rm(\cite{Z})}\label{multi-set}
For an integer $r\geq2$, let $K_{n_{1},n_{2},...,n_{r}}$ be the complete $r$-partite graph of order $n$ with parts $V_{1},V_{2},...,V_{r}$, where $n_{i}=|V_{i}|\geq1$ for $1\leq i\leq r$ and $\sum_{1\leq i\leq r}n_{i}=n$. For each $1\leq i\leq r$, let $H_{i}$ be a graph with vertex set $V(H_{i})\subseteq V_{i}$. Let $G$ be the graph obtained from $K_{n_{1},n_{2},...,n_{r}}$ by embedding the edges of $H_{i}$ into $V_{i}$ for all $1\leq i\leq r$. Set $\rho=\rho(G)$. Then $$\sum_{1\leq s\leq r}\frac{1}{1+\frac{n_{s}}{\rho}+\sum^{\infty}_{i=1}\frac{W^{i}(H_{s})}{\rho^{i+1}}}=r-1.$$
\end{theorem}

The existence of the infinite series in Theorem \ref{multi-set} was guaranteed by the following two results. The first one is taken  from Theorem 8.1.1 of \cite{CRS}.

\begin{lem}{\rm (\cite{CRS})}\label{subgraph}
If $H$ is a subgraph of a connected graph $G$, then $\rho(H)\leq\rho(G)$, with equality if and only if $H=G$.
\end{lem}

\begin{lem}{\rm(\cite{Z})}\label{exist}
Let $G$ be a graph of order $n$.  Then $\sum^{\infty}_{k=1}\frac{W^{k}(G)}{x^{k}}$ exists for any $x>\rho(G)$.
\end{lem}

\section{Proofs of Lemmas \ref{1-path}, \ref{2-path} and \ref{3-path}}

Let $G$ be a graph.  For a vertex $u$ and a walk $Q$  of $G$, we say that $Q$ crosses $u$, if $u$ is one of the vertices in the sequence of $Q$. For any integer $\ell\geq1$ and $u,v\in V(G)$, let $W^{\ell}_{u,v}(G)$ be the number of walks $Q$ of length $\ell$ in $G$, such that $Q$ crosses $u$ and $v$.  Recall that $P_{0}$ denotes the empty graph.

\begin{lem}\label{w-difference}
For integers $n_{1}\geq n_{2}+2\geq3$, denote $P_{n_{1}}=u_{1}u_{2}\cdots u_{n_{1}}$. Then $W^{\ell}(P_{n_{1}}\cup P_{n_{2}})-W^{\ell}(P_{n_{1}-1}\cup P_{n_{2}+1})
=W^{\ell}_{u_{1},u_{n_{2}+2}}(P_{n_{1}})$ for any $\ell\geq1$.
\end{lem}

\f{\bf Proof:} Let $G_{1}=P_{n_{1}}\cup P_{n_{2}}$ and $G_{2}=P_{n_{1}-1}\cup P_{n_{2}+1}$.
For any given $\ell\geq1$, we shall prove 
$$W^{\ell}(G_{1})-W^{\ell}(G_{2})
=W^{\ell}_{u_{1},u_{n_{2}+2}}(P_{n_{1}}).$$
 Recall that $P_{n_{1}}=u_{1}u_{2}\cdots u_{n_{1}}$ in $G_{1}$. In $G_{2}$, denote $P_{n_{2}+1}=v_{1}v_{2}\cdots v_{n_{2}}v_{n_{2}+1}$. Let $X$ denote the set of walks $Q$ of length $\ell$ in $G_{1}$, such that $Q$  crosses $u_{1}$. Let $Y$ denote the set of walks $Q$ of length $\ell$ in $G_{2}$, such that $Q$ crosses $v_{1}$. Clearly, $W^{\ell}(G_{1})=|X|+W^{\ell}(P_{n_{1}-1}\cup P_{n_{2}})$ and $W^{\ell}(G_{2})=|Y|+W^{\ell}(P_{n_{1}-1}\cup P_{n_{2}})$. 
This means that $$W^{\ell}(G_{1})-W^{\ell}(G_{2})=|X|-|Y|.$$
Now we divide $X$ in two parts: let $X_{1}$ denote the set of walks $Q$ of length $\ell$ in $G_{1}$, such that $Q$  crosses $u_{1}$ and $u_{n_{2}+2}$; let $X_{2}$ denote the set of walks $Q$ of length $\ell$ in $G_{1}$, such that $Q$  crosses $u_{1}$ and does not cross $u_{n_{2}+2}$. 
Clearly, $X=X_{1}\cup X_{2}$. By definitions of the sets $X_{1},X_{2},Y$, we see $|X_{2}|=|Y|$ and $|X_{1}|=W^{\ell}_{u_{1},u_{n_{2}+2}}(P_{n_{1}})$. Thus $|X|-|Y|=|X_{1}|=W^{\ell}_{u_{1},u_{n_{2}+2}}(P_{n_{1}})$. Hence, $W^{\ell}(G_{1})-W^{\ell}(G_{2})
=W^{\ell}_{u_{1},u_{n_{2}+2}}(P_{n_{1}})$. This completes the proof. \hfill$\Box$

\medskip

The following result will be used.

\medskip

\begin{lem}\label{sum}
 For integers $n\geq6m\geq0$, we have
 $$\binom{n}{0}4^{m}+\binom{n}{1}4^{m-1}+\cdots+\binom{n}{m}4^{0}\leq\frac{5}{2}\frac{n^{m}}{m!}$$
\end{lem}

\f{\bf Proof:} To prove the lemma, it is equivalent to show 
$$\binom{n}{0}+\binom{n}{1}\frac{1}{4}+\binom{n}{2}\frac{1}{4^{2}}+\cdots+\binom{n}{m}\frac{1}{4^{m}}
\leq\frac{5}{2}\frac{n^{m}}{4^{m}m!}.$$
We will prove this inequality by induction on $m$. It is easy to check for $m=0,1$. When $m=2$,
we need to show $\binom{n}{0}+\binom{n}{1}\frac{1}{4^{1}}+\binom{n}{2}\frac{1}{4^{2}}\leq\frac{5}{2}\frac{n^{2}}{4^{2}2}$, where $n\geq12$. This is equivalent to show $\frac{3}{2}n^{2}-3n-16\geq0$, which is clearly true as $n\geq12$.
When $m=3$,
we need to show $\binom{n}{0}+\binom{n}{1}\frac{1}{4^{1}}+\binom{n}{2}\frac{1}{4^{2}}+\binom{n}{3}\frac{1}{4^{3}}
\leq\frac{5}{2}\frac{n^{3}}{4^{3}6}$, where $n\geq18$. This is equivalent to show $3n^{3}-18n^{2}-172n-768\geq0$, which is clearly true as $n\geq18$.

Now we suppose the inequality  holds for $m=k$, where $k\geq3$. We will prove the inequality for $m=k+1$. Now let $m=k+1$, implying $n\geq6(k+1)$. Since the inequality  holds for $m=k$, we have  
$$\binom{n}{0}+\binom{n}{1}\frac{1}{4}+\binom{n}{2}\frac{1}{4^{2}}+\cdots+\binom{n}{k}\frac{1}{4^{k}}
\leq\frac{5}{2}\frac{n^{k}}{4^{k}k!}.$$
Then $$\binom{n}{0}+\binom{n}{1}\frac{1}{4}+\binom{n}{2}\frac{1}{4^{2}}+\cdots+\binom{n}{k}\frac{1}{4^{k}}
+\binom{n}{k+1}\frac{1}{4^{k+1}}\leq\binom{n}{k+1}\frac{1}{4^{k+1}}+\frac{5}{2}\frac{n^{k}}{4^{k}k!}.$$
To prove the inequality for $m=k+1$, it suffices to show 
$$\binom{n}{k+1}\frac{1}{4^{k+1}}+\frac{5}{2}\frac{n^{k}}{4^{k}k!}\leq
\frac{5}{2}\frac{n^{k+1}}{4^{k+1}(k+1)!}.$$
This is equivalent to show 
$$(n-1)(n-2)\cdots(n-k)\leq
\frac{5}{2}n^{k-1}(n-4k-4).$$ 
It suffices to show 
$$(n-k+1)(n-k)\leq
\frac{5}{2}n(n-4k-4).$$ 
This is equivalent to show
$$n(3n-16k-22)\geq2k^{2}-2k,$$
which is clearly true as $n\geq6k+6$ and $k\geq3$.
Thus, the inequality also holds for $m=k+1$.
By induction on $m$, this completes the proof. \hfill$\Box$

\begin{lem}\label{w-evaluation}
Let $G$ be a path of order $n\geq3$. Let $u,v$ be two vertices of $G$ such that $d_{G}(u,v)=\ell$, where $2\leq\ell\leq n-1$. Then $$\sum_{a\geq1}\frac{W^{a}_{u,v}(G)}{x^{a}}\leq\frac{1}{x^{\ell}}(\frac{30x}{x-1}e^{\frac{\lfloor\frac{3\ell}{2}\rfloor}{x^{2}}}
+\frac{32x}{x-8})$$ 
for any $x\geq\max\left\{\sqrt{n},9\right\}$.
\end{lem}

\begin{figure}
  \centering
  \includegraphics[width=11cm]{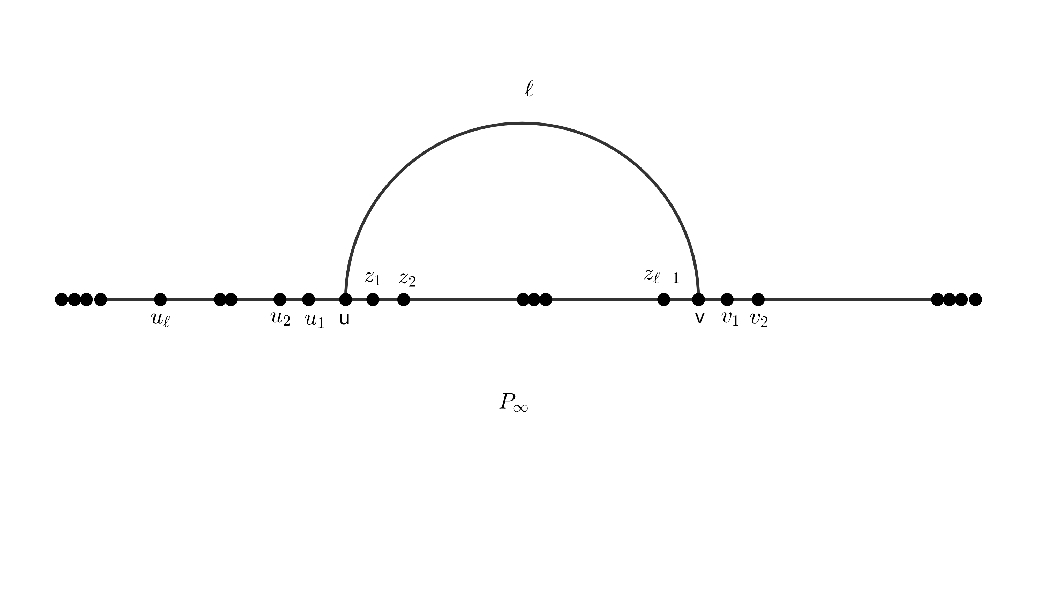}
  \caption{The infinite path $P_{\infty}$}\label{figure1}
\end{figure}

\f{\bf Proof:} Let $P_{\infty}$ be a path of "infinite" order, such that $u,v$ are two vertices of distance $\ell$ in $P_{\infty}$. That is, 
$$P_{\infty}=\cdots u_{\ell}\cdots u_{2}u_{1}uz_{1}z_{2}\cdots z_{\ell-1}vv_{1}v_{2}\cdots.$$ 
(See Figure \ref{figure1}, where the "semicircle" denotes the distance $\ell$ between $u$ and $v$.) Clearly, $W^{s}_{u,v}(G)\leq W^{s}_{u,v}(P_{\infty})$ for any $s\geq1$. For any vertex $y\in P_{\infty}$ and integer $a\geq1$, let $w^{a}_{uv}(y)$ be the number of walks $Q$ of length $a$ in $P_{\infty}$, such that $Q$ starts at $y$ and crosses $u$ and $v$.  

\medskip

\f{\bf Claim 1.} Let $1\leq a\leq\frac{3\ell}{2}$ be an integer. Then $w^{a}_{uv}(z_{i})\leq w^{a}_{uv}(u_{i})$ and $w^{a}_{uv}(z_{\ell-i})\leq w^{a}_{uv}(v_{i})$ for any $1\leq i< \frac{\ell}{2}$. Moreover, $w^{a}_{uv}(z_{\frac{\ell}{2}})\leq w^{a}_{uv}(u_{\frac{\ell}{2}})+w^{a}_{uv}(v_{\frac{\ell}{2}})$ for even $\ell$.

\medskip

\f{\bf Proof of Claim 1.} Clearly, $w^{a}_{uv}(z_{i})=0$ for $1\leq a\leq \ell-1$. Now assume $\ell\leq a\leq\frac{3\ell}{2}$. Let $1\leq i_{0}< \frac{\ell}{2}$ be fixed. We will show that $$w^{a}_{uv}(z_{i_{0}})\leq w^{a}_{uv}(u_{i_{0}}).$$
Let $X$ be the set of walks $Q$ starting at $z_{i_{0}}$ of length $a$ in $P_{\infty}$, such that $Q$ crosses $u$ and $v$.   Clearly, $w^{a}_{uv}(z_{i_{0}})=|X|$. Note that $1\leq i_{0}< \frac{\ell}{2}$ and $ a\leq\frac{3\ell}{2}$. Each $Q\in X$ crosses $u$ first.

Now we show that $|X|\leq w^{a}_{uv}(u_{i_{0}})$. Let $D$ be the set of walks $Q$ starting at $u_{i_{0}}$ of length $a$ in $P_{\infty}$, such that $Q$ crosses $u$ and $v$. Clearly, $w^{a}_{uv}(u_{i_{0}})=|D|$. We can define an injective map $\phi$ between $X$ and $D$. 
For any $Q\in X$, since $Q$ crosses $u$ first, we can denote 
$$Q=z_{i_{0}}z_{j_{1}}\cdots z_{j_{s}}uf_{1}f_{2}\cdots,$$
 where  $1\leq j_{1},...,j_{s}\leq \ell-1$ and $f_{1},f_{2},...$ are vertices of $P_{\infty}$. 
Define $$\phi(Q)=u_{i_{0}}u_{j_{1}}\cdots u_{j_{s}}uf_{1}f_{2}\cdots.$$
 It is easy to see that $\phi$ is injective. Hence, $|X|\leq|D|=w^{a}_{uv}(u_{i_{0}})$. That is, $w^{a}_{uv}(z_{i_{0}})\leq w^{a}_{uv}(u_{i_{0}})$. Similarly, we can show that $w^{a}_{uv}(z_{\ell-i_{0}})\leq w^{a}_{uv}(v_{i_{0}})$. 

Now assume that $\ell\geq2$ is even. Let $Y_{u}$ be the set of walks $Q$ starting at $z_{\frac{\ell}{2}}$ of length $a$ in $P_{\infty}$, such that $Q$ crosses $u$ first before $v$.  Let $Y_{v}$ be the set of walks $Q$ starting at $z_{\frac{\ell}{2}}$ of length $a$ in $P_{\infty}$, such that $Q$ crosses $v$ first before $u$. Clearly, $w^{a}_{uv}(z_{\frac{\ell}{2}})=|Y_{u}|+|Y_{v}|$. Similar to the above discussion, we can show that $|Y_{u}|\leq w^{a}_{uv}(u_{\frac{\ell}{2}})$ and $|Y_{v}|\leq w^{a}_{uv}(v_{\frac{\ell}{2}})$. Hence $w^{a}_{uv}(z_{\frac{\ell}{2}})\leq w^{a}_{uv}(u_{\frac{\ell}{2}})+w^{a}_{uv}(v_{\frac{\ell}{2}})$.
 This finishes the proof of Claim 1. \hfill$\Box$

\medskip

By Claim 1, we see 
$$\sum_{1\leq i\leq \ell-1}w^{a}_{uv}(z_{i})\leq\sum_{1\leq i\leq \frac{\ell}{2}}w^{a}_{uv}(u_{i})+\sum_{1\leq i\leq \frac{\ell}{2}}w^{a}_{uv}(v_{i})$$
 for any $1\leq a\leq\frac{3}{2}\ell$. Clearly, $w^{a}_{uv}(u_{i})=w^{a}_{uv}(v_{i})=0$ for any $i>\frac{\ell}{2}$ and $a\leq\frac{3}{2}\ell$. 
 Thus, for any $a\leq\frac{3}{2}\ell$, 
 \begin{equation}
\begin{aligned}
W^{a}_{u,v}(P_{\infty})&=w^{a}_{uv}(u)+w^{a}_{uv}(v)+\sum_{1\leq i\leq \ell-1}w^{a}_{uv}(z_{i})+\sum_{1\leq i\leq \frac{\ell}{2}}w^{a}_{uv}(u_{i})+\sum_{1\leq i\leq \frac{\ell}{2}}w^{a}_{uv}(v_{i})\\
&\leq 2\left(w^{a}_{uv}(u)+w^{a}_{uv}(v)+\sum_{1\leq i\leq \frac{\ell}{2}}w^{a}_{uv}(u_{i})+\sum_{1\leq i\leq \frac{\ell}{2}}w^{a}_{uv}(v_{i})\right).
\end{aligned}\notag
\end{equation}
 By symmetry, $w^{a}_{uv}(u)=w^{a}_{uv}(v)$ and $w^{a}_{uv}(u_{i})=w^{a}_{uv}(v_{i})$ for any $a\geq1$ and $i\geq1$. Thus, for any $a\leq\frac{3\ell}{2}$, 
 $$W^{a}_{u,v}(P_{\infty})\leq4(w^{a}_{uv}(u)+\sum_{1\leq i\leq \frac{\ell}{2}}w^{a}_{uv}(u_{i})).$$
 
 \medskip

\f{\bf Claim 2.} Let $a,i$ be given integers with $\ell\leq a\leq\frac{3\ell}{2}$ and $1\leq i\leq \frac{\ell}{2}$. Then
$$w^{a}_{uv}(u_{i})\leq\frac{5}{2}\frac{(\lfloor\frac{3\ell}{2}\rfloor)^{\lfloor\frac{a-\ell-i}{2}\rfloor}}{\lfloor\frac{a-\ell-i}{2}\rfloor !}$$ for even $a-\ell-i$,
and
$$w^{a}_{uv}(u_{i})\leq\frac{5(\lfloor\frac{3\ell}{2}\rfloor)^{\lfloor\frac{a-\ell-i}{2}\rfloor}}{\lfloor\frac{a-\ell-i}{2}\rfloor !}$$ for odd $a-\ell-i$.

\medskip

\f{\bf Proof of Claim 2.} Since $P_{\infty}$ is a path, each walk starting at $u_{i}$ corresponds to a "$\pm$"-series: $$u_{i}+--+-+\cdots.$$
  (For example, $u_{1}+-++--$ corresponds to the walk $u_{1}uu_{1}uz_{1}uu_{1}$.) Recall that $\ell\leq a\leq\frac{3\ell}{2}$ and $1\leq i\leq \frac{\ell}{2}$. Note that the distance between $u_{i}$ and $v$ is $\ell+i$. For any $\ell+i\leq j\leq a$, let $X_{j}$ be the set of walks $Q$ starting at $u_{i}$ of length $a$ in $P_{\infty}$, such that $Q$ crosses $v$ at step $j$ (with smallest value). That is, $Q=u_{i}g_{1}g_{2}\cdots g_{a}$, where $g_{j}=v$ and $g_{k}\neq v$ for any $1\leq k\leq j-1$.  Clearly, $j\equiv \ell+i(\mod~2)$. Moreover, $$w^{a}_{uv}(u_{i})=|X_{\ell+i}|+|X_{\ell+i+2}|+\cdots+|X_{a}|$$
   for even $a-\ell-i$, and 
   $$w^{a}_{uv}(u_{i})=|X_{\ell+i}|+|X_{\ell+i+2}|+\cdots+|X_{a-1}|$$
    for odd $a-\ell-i$.
  
 For any $\ell+i\leq j\leq a$ and $j\equiv \ell+i(\mod~2)$, in the corresponding "$\pm$"-series to walks $Q\in X_{j}$, there are exactly $\lfloor\frac{j-\ell-i}{2}\rfloor$ symbols "$-$" in the first $j$ symbols "$\pm$". Hence, 
 $$|X_{j}|\leq\binom{j}{\lfloor\frac{j-\ell-i}{2}\rfloor}2^{a-j}.$$
  Since $j\leq a\leq\lfloor\frac{3\ell}{2}\rfloor$, we have $$|X_{j}|\leq\binom{\lfloor\frac{3\ell}{2}\rfloor}{\lfloor\frac{j-\ell-i}{2}\rfloor}2^{a-j}.$$ Let $m=\lfloor\frac{a-\ell-i}{2}\rfloor$. Clearly, $m\leq\frac{a-\ell-i}{2}\leq\frac{\ell}{4}$ as $a\leq\frac{3\ell}{2}$. Then $6m\leq\frac{3\ell}{2}$, implying $6m\leq\lfloor\frac{3\ell}{2}\rfloor$.
  
  For even $a-\ell-i$,
 using Lemma \ref{sum} as $6m\leq\lfloor\frac{3\ell}{2}\rfloor$,
  \begin{equation}
\begin{aligned}
w^{a}_{uv}(u_{i})&=|X_{\ell+i}|+|X_{\ell+i+2}|+\cdots+|X_{a}|\\
&\leq\binom{\lfloor\frac{3\ell}{2}\rfloor}{0}2^{a-\ell-i}+
  \binom{\lfloor\frac{3\ell}{2}\rfloor}{1}2^{a-\ell-i-2}+ \binom{\lfloor\frac{3\ell}{2}\rfloor}{2}2^{a-\ell-i-4}+
  \cdots+
  \binom{\lfloor\frac{3\ell}{2}\rfloor}{\frac{a-\ell-i}{2}}2^{0}\\
&=\binom{\lfloor\frac{3\ell}{2}\rfloor}{0}4^{m}+
  \binom{\lfloor\frac{3\ell}{2}\rfloor}{1}4^{m-1}+ 
  \cdots+
  \binom{\lfloor\frac{3\ell}{2}\rfloor}{m}4^{0}\\
&\leq\frac{5}{2}\frac{(\lfloor\frac{3\ell}{2}\rfloor)^{m}}{m !}\\
&=
\frac{5}{2}\frac{(\lfloor\frac{3\ell}{2}\rfloor)^{\lfloor\frac{a-\ell-i}{2}\rfloor}}{\lfloor\frac{a-\ell-i}{2}\rfloor !}.
\end{aligned}\notag
\end{equation}
For odd $a-\ell-i$, we have $m=\frac{a-1-\ell-i}{2}$.
 Using Lemma \ref{sum} as $6m\leq\lfloor\frac{3\ell}{2}\rfloor$,
  \begin{equation}
\begin{aligned}
w^{a}_{uv}(u_{i})&=|X_{\ell+i}|+|X_{\ell+i+2}|+\cdots+|X_{a-1}|\\
&\leq\binom{\lfloor\frac{3\ell}{2}\rfloor}{0}2^{a-\ell-i}+
  \binom{\lfloor\frac{3\ell}{2}\rfloor}{1}2^{a-\ell-i-2}+ \binom{\lfloor\frac{3\ell}{2}\rfloor}{2}2^{a-\ell-i-4}+
  \cdots+
  \binom{\lfloor\frac{3\ell}{2}\rfloor}{\lfloor\frac{a-\ell-i}{2}\rfloor}2^{1}\\
 &= 2\binom{\lfloor\frac{3\ell}{2}\rfloor}{0}2^{a-1-\ell-i}+
  2\binom{\lfloor\frac{3\ell}{2}\rfloor}{1}2^{a-1-\ell-i-2}+ 2\binom{\lfloor\frac{3\ell}{2}\rfloor}{2}2^{a-1-\ell-i-4}+
  \cdots+
  2\binom{\lfloor\frac{3\ell}{2}\rfloor}{\lfloor\frac{a-\ell-i}{2}\rfloor}2^{0}\\
&=2\binom{\lfloor\frac{3\ell}{2}\rfloor}{0}4^{m}+
  2\binom{\lfloor\frac{3\ell}{2}\rfloor}{1}4^{m-1}+ 
  \cdots+
  2\binom{\lfloor\frac{3\ell}{2}\rfloor}{m}4^{0}\\
&\leq\frac{5(\lfloor\frac{3\ell}{2}\rfloor)^{m}}{m !}\\
&=
\frac{5(\lfloor\frac{3\ell}{2}\rfloor)^{\lfloor\frac{a-\ell-i}{2}\rfloor}}{\lfloor\frac{a-\ell-i}{2}\rfloor !}.
\end{aligned}\notag
\end{equation}
This finishes the proof of Claim 2. \hfill$\Box$

\medskip

Similar to Claim 2, we can show $w^{a}_{uv}(u)\leq\frac{5}{2}\frac{(\lfloor\frac{3\ell}{2}\rfloor)^{\lfloor\frac{a-\ell}{2}\rfloor}}{\lfloor\frac{a-\ell}{2}\rfloor !}$ for even $a-\ell$, and $w^{a}_{uv}(u)\leq\frac{5(\lfloor\frac{3\ell}{2}\rfloor)^{\lfloor\frac{a-\ell}{2}\rfloor}}{\lfloor\frac{a-\ell}{2}\rfloor !}$ for odd $a-\ell$.
Clearly, $w^{a}_{uv}(u_{i})=0$ for $i>a-\ell$. By Claim 2, for any $\ell\leq a\leq \frac{3\ell}{2}$,
\begin{equation}
\begin{aligned}
&w^{a}_{uv}(u)+\sum_{1\leq i\leq\frac{\ell}{2}}w^{a}_{uv}(u_{i})\\
&=w^{a}_{uv}(u)+\sum_{1\leq i\leq a-\ell}w^{a}_{uv}(u_{i})\\
 &\leq(5+\frac{5}{2})\times(1+\lfloor\frac{3\ell}{2}\rfloor+
\frac{(\lfloor\frac{3\ell}{2}\rfloor)^{2}}{2!}+\cdots+
\frac{(\lfloor\frac{3\ell}{2}\rfloor)^{\lfloor\frac{a-\ell}{2}\rfloor}}{\lfloor\frac{a-\ell}{2}\rfloor!})\\
&=\frac{15}{2}\sum_{0\leq i\leq\frac{a-\ell}{2}}\frac{(\lfloor\frac{3\ell}{2}\rfloor)^{i}}{i!}.
\end{aligned}\notag
\end{equation}
Then
\begin{equation}
\begin{aligned}
\sum_{\ell\leq a\leq\frac{3\ell}{2}}\frac{W^{a}_{u,v}(P_{\infty})}{x^{a}}&\leq\sum_{\ell\leq a\leq\frac{3\ell}{2}}\frac{4}{x^{a}}(w^{a}_{uv}(u)+\sum_{1\leq i\leq \frac{\ell}{2}}w^{a}_{uv}(u_{i}))\\
&\leq\frac{30}{x^{\ell}}\sum_{\ell\leq a\leq\frac{3\ell}{2}}\frac{1}{x^{a-\ell}}\sum_{0\leq i\leq \frac{a-\ell}{2}}\frac{(\lfloor\frac{3\ell}{2}\rfloor)^{i}}{i!}\\
&=\frac{30}{x^{\ell}}\sum_{0\leq i\leq\frac{\ell}{4}}\frac{(\lfloor\frac{3\ell}{2}\rfloor)^{i}}{i!}\sum_{\ell+2i\leq a\leq \frac{3\ell}{2}}\frac{1}{x^{a-\ell}}\\
&\leq\frac{30}{x^{\ell}}\sum_{0\leq i\leq\frac{\ell}{4}}\frac{(\lfloor\frac{3\ell}{2}\rfloor)^{i}}{i!}\frac{1}{x^{2i}}\frac{x}{x-1}\\
&\leq\frac{30}{x^{\ell}}\frac{x}{x-1}\sum^{\infty}_{i=0}\frac{1}{i!}(\frac{\lfloor\frac{3\ell}{2}\rfloor}{x^{2}})^{i}\\
&=\frac{30}{x^{\ell}}\frac{x}{x-1}e^{\frac{\lfloor\frac{3\ell}{2}\rfloor}{x^{2}}}.
\end{aligned}\notag
\end{equation}
Note that $W^{s}_{u,v}(G)=0$ for any $s\leq \ell-1$ as $d_{G}(u,v)=\ell$. Then
$$\sum_{1\leq a\leq\frac{3\ell}{2}}\frac{W^{a}_{u,v}(G)}{x^{a}}\leq\sum_{\ell\leq a\leq\frac{3\ell}{2}}\frac{W^{a}_{u,v}(P_{\infty})}{x^{a}}\leq
\frac{30}{x^{\ell}}\frac{x}{x-1}e^{\frac{\lfloor\frac{3\ell}{2}\rfloor}{x^{2}}}.$$

Since $G$ has maximum degree 2, we see $W^{a}_{u,v}(G)\leq2^{a}n$ for any $a\geq1$.
Note that $a-5\leq3(a-\ell-2)$ and $a-\ell-2\geq0$, if $a>\frac{3\ell}{2}$ and $\ell\geq2$. 
Thus,
\begin{equation}
\begin{aligned}
\sum_{a>\frac{3\ell}{2}}\frac{W^{a}_{u,v}(G)}{x^{a}}&\leq\sum_{a>\frac{3\ell}{2}}\frac{2^{a}n}{x^{a}}\\
&\leq\frac{32}{x^{\ell}}\sum_{a>\frac{3\ell}{2}}\frac{2^{a-5}}{x^{a-\ell-2}}~~~(using~x\geq\sqrt{n})\\
&\leq\frac{32}{x^{\ell}}\sum_{a>\frac{3\ell}{2}}(\frac{8}{x})^{a-\ell-2}~~~(using~a-5\leq3(a-\ell-2))\\
&\leq\frac{32}{x^{\ell}}\sum_{j\geq0}(\frac{8}{x})^{j}~~~(using~a-\ell-2\geq0)\\
&=\frac{32}{x^{\ell}}\frac{x}{x-8}~~~(using~x\geq9).
\end{aligned}\notag
\end{equation}
Consequently, $$\sum_{a\geq1}\frac{W^{a}_{u,v}(G)}{x^{a}}\leq\frac{1}{x^{\ell}}(\frac{30x}{x-1}e^{\frac{\lfloor\frac{3\ell}{2}\rfloor}{x^{2}}}
+\frac{32x}{x-8}).$$
This completes the proof. \hfill$\Box$

\medskip

Now we are ready to prove the three lemmas.

\medskip

\f{\bf Proof of Lemma \ref{1-path}.}  Let $G^{1}=P_{n_{1}}\cup P_{n_{2}}\cup T$ and $G^{2}=P_{n_{1}-1}\cup P_{n_{2}+1}\cup T$. Then $G_{j}=H\vee G^{j}$ for $j=1,2$. Set $\rho(G_{j})=\rho_{j}$ for any $1\leq j\leq2$. We will prove $\rho_{1}>\rho_{2}$. 
By Theorem \ref{multi-set}, $$\frac{1}{1+\frac{|H|}{\rho_{j}}+\frac{1}{\rho_{j}}\sum^{\infty}_{i=1}\frac{W^{i}(H)}{\rho^{i}_{j}}}
+\frac{1}{1+\frac{|G^{j}|}{\rho_{j}}+\frac{1}{\rho_{j}}\sum^{\infty}_{i=1}\frac{W^{i}(G^{j})}{\rho^{i}_{j}}}=1,$$
where $1\leq j\leq2$. 
For $1\leq j\leq2$, let 
$$f(G_{j},x)=\frac{1}{1+\frac{|H|}{x}+\frac{1}{x}\sum^{\infty}_{i=1}\frac{W^{i}(H)}{x^{i}}}
+\frac{1}{1+\frac{|G^{j}|}{x}+\frac{1}{x}\sum^{\infty}_{i=1}\frac{W^{i}(G^{j})}{x^{i}}}.$$
 (Clearly, the domain of $f(G_{j},x)$ is of the form $x\geq a_{j}$ for some positive number $a_{j}$.)
  Clearly, $f(G_{j},x)$ is strictly increasing with respective to $x$ in its domain. Thus $\rho_{j}$ is the largest root of $f(G_{j},x)=1$.

In $G^{1}$, denote $P_{n_{1}}=u_{1}u_{2}\cdots u_{n_{1}}$. For any $\ell\geq1$, by Lemma \ref{w-difference}, 
$$W^{\ell}(P_{n_{1}}\cup P_{n_{2}})-W^{\ell}(P_{n_{1}-1}\cup P_{n_{2}+1})
=W^{\ell}_{u_{1},u_{n_{2}+2}}(P_{n_{1}}).$$
 This means that $W^{\ell}(G^{1})-W^{\ell}(G^{2})=W^{\ell}_{u_{1},u_{n_{2}+2}}(P_{n_{1}})\geq0$. Clearly, $W^{n_{2}+1}_{u_{1},u_{n_{2}+2}}(P_{n_{1}})=2$. It follows that $$\sum^{\infty}_{i=1}\frac{W^{i}(G^{1})}{\rho_{2}^{i}}>
\sum^{\infty}_{i=1}\frac{W^{i}(G^{2})}{\rho_{2}^{i}}.$$
 This implies that $$f(G_{1},\rho_{2})<f(G_{2},\rho_{2}).$$
  Note that $f(G_{2},\rho_{2})=1$.
 So, $f(G_{1},\rho_{2})<1$. Recall that $f(G_{1},x)$ is strictly increasing in its domain, and $f(G_{1},\rho_{1})=1$. We must have $\rho_{1}>\rho_{2}$. This completes the proof. \hfill$\Box$

\medskip

\f{\bf Proof of Lemma \ref{2-path}.}  Let $G^{1}=P_{n_{1}}\cup P_{n_{2}}\cup P_{n_{3}}\cup P_{n_{4}}\cup T$ and $G^{2}=P_{n_{1}-1}\cup P_{n_{2}+1}\cup P_{n_{3}+1}\cup P_{n_{4}-1}\cup T$. Then $G_{j}=H\vee G^{j}$ for $j=1,2$. Set $\rho(G_{j})=\rho_{j}$ for any $1\leq j\leq2$. We will prove $\rho_{1}<\rho_{2}$. Since $G_{j}$ contains $K_{1}\vee\overline{K_{|T|+\sum_{1\leq i\leq4}n_{i}}}$ as a subgraph, we have 
$$\rho_{j}\geq\rho(K_{1}\vee\overline{K_{|T|+\sum_{1\leq i\leq4}n_{i}}})=\sqrt{|T|+\sum_{1\leq i\leq4}n_{i}}\geq\max\left\{\sqrt{2000},\sqrt{\sum_{1\leq i\leq4}n_{i}}\right\}.$$ 
By Theorem \ref{multi-set}, $$\frac{1}{1+\frac{|H|}{\rho_{j}}+\frac{1}{\rho_{j}}\sum^{\infty}_{i=1}\frac{W^{i}(H)}{\rho^{i}_{j}}}
+\frac{1}{1+\frac{|G^{j}|}{\rho_{j}}+\frac{1}{\rho_{j}}\sum^{\infty}_{i=1}\frac{W^{i}(G^{j})}{\rho^{i}_{j}}}=1,$$
where $1\leq j\leq2$. 
For $1\leq j\leq2$, let $$f(G_{j},x)=\frac{1}{1+\frac{|H|}{x}+\frac{1}{x}\sum^{\infty}_{i=1}\frac{W^{i}(H)}{x^{i}}}
+\frac{1}{1+\frac{|G^{j}|}{x}+\frac{1}{x}\sum^{\infty}_{i=1}\frac{W^{i}(G^{j})}{x^{i}}}.$$
  Similar to Lemma \ref{1-path}, $f(G_{j},x)$ is strictly increasing with respective to $x$ in its domain, and $f(G_{j},\rho_{j})=1$.

In $G^{1}$, denote $P_{n_{1}}=u_{1}u_{2}\cdots u_{n_{1}}$. In $G^{2}$, denote $P_{n_{3}+1}=v_{1}v_{2}\cdots v_{n_{3}+1}$. For any $\ell\geq1$, by Lemma \ref{w-difference}, $$W^{\ell}(P_{n_{1}}\cup P_{n_{2}})-W^{\ell}(P_{n_{1}-1}\cup P_{n_{2}+1})
=W^{\ell}_{u_{1},u_{n_{2}+2}}(P_{n_{1}}),$$ and 
$$W^{\ell}(P_{n_{3}+1}\cup P_{n_{4}-1})-W^{\ell}(P_{n_{3}}\cup P_{n_{4}})
=W^{\ell}_{v_{1},v_{n_{4}+1}}(P_{n_{3}+1}).$$
 This means that 
$$W^{\ell}(G^{1})-W^{\ell}(G^{2})=
W^{\ell}_{u_{1},u_{n_{2}+2}}(P_{n_{1}})-W^{\ell}_{v_{1},v_{n_{4}+1}}(P_{n_{3}+1}).$$
Note that $W^{i}_{u_{1},u_{n_{2}+2}}(P_{n_{1}})\leq W^{i}_{u_{1},u_{n_{4}+2}}(P_{n_{1}})$ for any $i\geq1$ by the definition of $W^{i}_{u_{1},u_{n_{2}+2}}(P_{n_{1}})$, since $n_{2}\geq n_{4}$.
Thus, $$W^{\ell}(G^{1})-W^{\ell}(G^{2})\leq
W^{\ell}_{u_{1},u_{n_{4}+2}}(P_{n_{1}})-W^{\ell}_{v_{1},v_{n_{4}+1}}(P_{n_{3}+1}).$$
Then 
$$\sum^{\infty}_{i=1}(\frac{W^{i}(G^{1})}{\rho_{1}^{i}}-\frac{W^{i}(G^{2})}{\rho_{1}^{i}})=
\sum^{\infty}_{i=1}\frac{W^{i}(G^{1})-W^{i}(G^{2})}{\rho_{1}^{i}}\leq
\sum^{\infty}_{i=1}\frac{
W^{i}_{u_{1},u_{n_{4}+2}}(P_{n_{1}})-W^{i}_{v_{1},v_{n_{4}+1}}(P_{n_{3}+1})}{\rho_{1}^{i}}.$$

 Clearly, $W^{n_{4}}_{v_{1},v_{n_{4}+1}}(P_{n_{3}+1})=2$ and $W^{n_{4}+1}_{v_{1},v_{n_{4}+1}}(P_{n_{3}+1})\geq2$. 
 It follows that $$\sum^{\infty}_{i=1}\frac{W^{i}_{v_{1},v_{n_{4}+1}}(P_{n_{3}+1})}{\rho_{1}^{i}}\geq
 \frac{2}{\rho_{1}^{n_{4}}}+\frac{2}{\rho_{1}^{n_{4}+1}}=
 \frac{1}{\rho_{1}^{n_{4}+1}}(2\rho_{1}+2).$$
 By Lemma \ref{w-evaluation} (by letting $\ell=n_{4}+1$ there),
 $$\sum^{\infty}_{i=1}\frac{W^{i}_{u_{1},u_{n_{4}+2}}(P_{n_{1}})}{\rho_{1}^{i}}
 \leq\frac{1}{\rho_{1}^{n_{4}+1}}(\frac{30\rho_{1}}{\rho_{1}-1}e^{\frac{\frac{3}{2}(n_{4}+1)}{\rho_{1}^{2}}}+
\frac{32\rho_{1}}{\rho_{1}-8}).$$
Since $\rho_{1}\geq\max\left\{\sqrt{2000},\sqrt{\sum_{1\leq i\leq4}n_{i}}\right\}\geq\max\left\{\sqrt{2000},\sqrt{4n_{4}+2}\right\}$, we have $e^{\frac{\frac{3}{2}(n_{4}+1)}{\rho_{1}^{2}}}\leq e^{\frac{\frac{3}{2}(n_{4}+1)}{4n_{4}+2}}\leq e^{0.5}<1.65$.
Thus,
\begin{equation}
\begin{aligned}
&\sum^{\infty}_{i=1}(\frac{W^{i}(G^{1})}{\rho_{1}^{i}}-\frac{W^{i}(G^{2})}{\rho_{1}^{i}})\\
&\leq\frac{2}{\rho_{1}^{n_{4}+1}}(\frac{15\rho_{1}}{\rho_{1}-1}e^{\frac{\frac{3}{2}(n_{4}+1)}{\rho_{1}^{2}}}+
\frac{16\rho_{1}}{\rho_{1}-8}-\rho_{1}-1)\\
&\leq\frac{2}{\rho_{1}^{n_{4}+1}}(\frac{25\rho_{1}}{\rho_{1}-1}+
\frac{16\rho_{1}}{\rho_{1}-8}-\rho_{1}-1)~~~(using~e^{\frac{\frac{3}{2}(n_{4}+1)}{\rho_{1}^{2}}}<1.65)\\
&=\frac{2}{\rho_{1}^{n_{4}+1}}(40+\frac{25}{\rho_{1}-1}
+\frac{128}{\rho_{1}-8}-\rho_{1}).
\end{aligned}\notag
\end{equation}
 
Now consider $g(x)=40+\frac{25}{x-1}
+\frac{128}{x-8}-x$ for $x\geq\sqrt{2000}$. Clearly, $g(x)$ is decreasing for $x\geq\sqrt{2000}$. By a calculation, $g(\sqrt{2000})=-0.66\cdots<0$. Thus $g(x)<0$ for any $x\geq\sqrt{2000}$. Then $\sum^{\infty}_{i=1}(\frac{W^{i}(G^{1})}{\rho_{1}^{i}}-\frac{W^{i}(G^{2})}{\rho_{1}^{i}})<0$, or equivalently, $\sum^{\infty}_{i=1}\frac{W^{i}(G^{1})}{\rho_{1}^{i}}<
\sum^{\infty}_{i=1}\frac{W^{i}(G^{2})}{\rho_{1}^{i}}$. 
It follows that  $f(G_{1},\rho_{1})>f(G_{2},\rho_{1})$. Note that $f(G_{1},\rho_{1})=1$.
 So, $f(G_{2},\rho_{1})<1$. Recall that $f(G_{2},x)$ is strictly increasing in its domain, and $f(G_{2},\rho_{2})=1$. We must have $\rho_{1}<\rho_{2}$. This completes the proof. \hfill$\Box$

\medskip

\f{\bf Proof of Lemma \ref{3-path}.} Let $G^{1}=P_{n_{1}}\cup P_{n_{2}}\cup P_{n_{3}}\cup P_{n_{4}}\cup P_{n_{5}}\cup T$ and $G^{2}=P_{n_{1}-2}\cup P_{n_{2}+1}\cup P_{n_{3}+1}\cup P_{n_{4}+1}\cup P_{n_{5}-1}\cup T$. Then $G_{j}=H\vee G^{j}$ for $j=1,2$. Set $\rho(G_{j})=\rho_{j}$ for any $1\leq j\leq2$. We will prove $\rho_{1}<\rho_{2}$.
Since $G_{j}$ contains $K_{1}\vee\overline{K_{|T|+\sum_{1\leq i\leq5}n_{i}}}$ as a subgraph, we have 
$$\rho_{j}\geq\rho(K_{1}\vee\overline{K_{|T|+\sum_{1\leq i\leq5}n_{i}}})=\sqrt{|T|+\sum_{1\leq i\leq5}n_{i}}\geq\max\left\{\sqrt{6500},\sqrt{\sum_{1\leq i\leq5}n_{i}}\right\}.$$
 By Theorem \ref{multi-set}, $$\frac{1}{1+\frac{|H|}{\rho_{j}}+\frac{1}{\rho_{j}}\sum^{\infty}_{i=1}\frac{W^{i}(H)}{\rho^{i}_{j}}}
+\frac{1}{1+\frac{|G^{j}|}{\rho_{j}}+\frac{1}{\rho_{j}}\sum^{\infty}_{i=1}\frac{W^{i}(G^{j})}{\rho^{i}_{j}}}=1,$$
where $1\leq j\leq2$. 
For $1\leq j\leq2$, let $$f(G_{j},x)=\frac{1}{1+\frac{|H|}{x}+\frac{1}{x}\sum^{\infty}_{i=1}\frac{W^{i}(H)}{x^{i}}}
+\frac{1}{1+\frac{|G^{j}|}{x}+\frac{1}{x}\sum^{\infty}_{i=1}\frac{W^{i}(G^{j})}{x^{i}}}.$$
Clearly, $f(G_{j},x)$ is strictly increasing with respective to $x$ in its domain, and $f(G_{j},\rho_{j})=1$.

In $G^{1}$, denote $P_{n_{1}}=u_{1}u_{2}\cdots u_{n_{1}}$. In $G^{2}$, denote $P_{n_{4}+1}=v_{1}v_{2}\cdots v_{n_{4}+1}$. For any $\ell\geq1$, by Lemma \ref{w-difference}, $$W^{\ell}(P_{n_{1}}\cup P_{n_{2}}\cup P_{n_{3}})-W^{\ell}(P_{n_{1}-1}\cup P_{n_{2}+1}\cup P_{n_{3}})=W^{\ell}(P_{n_{1}}\cup P_{n_{2}})-W^{\ell}(P_{n_{1}-1}\cup P_{n_{2}+1})$$
$$=W^{\ell}_{u_{1},u_{n_{2}+2}}(P_{n_{1}})\leq W^{\ell}_{u_{1},u_{n_{5}+2}}(P_{n_{1}}).$$
 Similarly (denoting $P_{n_{1}-1}=u_{1}u_{2}\cdots u_{n_{1}-1}$),
$$W^{\ell}(P_{n_{1}-1}\cup P_{n_{2}+1}\cup P_{n_{3}})-W^{\ell}(P_{n_{1}-2}\cup P_{n_{2}+1}\cup P_{n_{3}+1})=W^{\ell}(P_{n_{1}-1}\cup P_{n_{3}})-W^{\ell}(P_{n_{1}-2}\cup P_{n_{3}+1})$$
$$= W^{\ell}_{u_{1},u_{n_{3}+2}}(P_{n_{1}-1})\leq W^{\ell}_{u_{1},u_{n_{5}+2}}(P_{n_{1}}).$$
It follows that 
$$W^{\ell}(P_{n_{1}}\cup P_{n_{2}}\cup P_{n_{3}})-W^{\ell}(P_{n_{1}-2}\cup P_{n_{2}+1}\cup P_{n_{3}+1})\leq2W^{\ell}_{u_{1},u_{n_{5}+2}}(P_{n_{1}}).$$
Also by Lemma \ref{w-difference}, 
$$W^{\ell}(P_{n_{4}+1}\cup P_{n_{5}-1})-W^{\ell}(P_{n_{4}}\cup P_{n_{5}})
=W^{\ell}_{v_{1},v_{n_{5}+1}}(P_{n_{4}+1}).$$ 
It follows that 
$$W^{\ell}(G^{1})-W^{\ell}(G^{2})\leq
2W^{\ell}_{u_{1},u_{n_{5}+2}}(P_{n_{1}})-W^{\ell}_{v_{1},v_{n_{5}+1}}(P_{n_{4}+1}).$$
Thus, 
$$\sum^{\infty}_{i=1}(\frac{W^{i}(G^{1})}{\rho_{1}^{i}}-\frac{W^{i}(G^{2})}{\rho_{1}^{i}})\leq
\sum^{\infty}_{i=1}\frac{2W^{\ell}_{u_{1},u_{n_{5}+2}}(P_{n_{1}})-
W^{\ell}_{v_{1},v_{n_{5}+1}}(P_{n_{4}+1})}{\rho_{1}^{i}}.$$

 Clearly, $W^{n_{5}}_{v_{1},v_{n_{5}+1}}(P_{n_{4}+1})=2$ and $W^{n_{5}+1}_{v_{1},v_{n_{5}+1}}(P_{n_{4}+1})\geq2$. 
 It follows that $$\sum^{\infty}_{i=1}\frac{W^{i}_{v_{1},v_{n_{5}+1}}(P_{n_{4}+1})}{\rho_{1}^{i}}\geq
 \frac{2}{\rho_{1}^{n_{5}}}+\frac{2}{\rho_{1}^{n_{5}+1}}=
 \frac{1}{\rho_{1}^{n_{5}+1}}(2\rho_{1}+2).$$ 
 By Lemma \ref{w-evaluation} (by letting $\ell=n_{5}+1$ there), 
   $$\sum^{\infty}_{i=1}\frac{2W^{i}_{u_{1},u_{n_{5}+2}}(P_{n_{1}})}{\rho_{1}^{i}}
 \leq\frac{2}{\rho_{1}^{n_{5}+1}}(\frac{30\rho_{1}}{\rho_{1}-1}e^{\frac{\frac{3}{2}(n_{5}+1)}{\rho_{1}^{2}}}+
\frac{32\rho_{1}}{\rho_{1}-8}).$$
Since $\rho_{1}\geq\max\left\{\sqrt{6500},\sqrt{\sum_{1\leq i\leq5}n_{i}}\right\}\geq\max\left\{\sqrt{6500},\sqrt{5n_{5}+3}\right\}$, we have $e^{\frac{\frac{3}{2}(n_{5}+1)}{\rho_{1}^{2}}}\leq e^{\frac{\frac{3}{2}(n_{5}+1)}{5n_{5}+3}}
<e^{0.4}<1.5$.
Thus,
\begin{equation}
\begin{aligned}
&\sum^{\infty}_{i=1}(\frac{W^{i}(G^{1})}{\rho_{1}^{i}}-\frac{W^{i}(G^{2})}{\rho_{1}^{i}})\\
&\leq\frac{2}{\rho_{1}^{n_{5}+1}}(\frac{30\rho_{1}}{\rho_{1}-1}e^{\frac{\frac{3}{2}(n_{5}+1)}{\rho_{1}^{2}}}+
\frac{32\rho_{1}}{\rho_{1}-8}-\rho_{1}-1)\\
&<\frac{2}{\rho_{1}^{n_{5}+1}}(\frac{45\rho_{1}}{\rho_{1}-1}+
\frac{32\rho_{1}}{\rho_{1}-8}-\rho_{1}-1)~~~(using~e^{\frac{2(n_{5}+1)}{\rho_{1}^{2}}}<1.5)\\
&=\frac{2}{\rho_{1}^{n_{5}+1}}(76+\frac{45}{\rho_{1}-1}
+\frac{256}{\rho_{1}-8}-\rho_{1}).
\end{aligned}\notag
\end{equation}
 
Now consider $g(x)=76+\frac{45}{x-1}
+\frac{256}{x-8}-x$ for $x\geq\sqrt{6500}$. Clearly, $g(x)$ is decreasing for $x\geq\sqrt{6500}$. By a calculation, $g(\sqrt{6500})=-0.53\cdots$. Thus $g(x)<0$ for any $x\geq\sqrt{6500}$. 
So, $\sum^{\infty}_{i=1}(\frac{W^{i}(G^{1})}{\rho_{1}^{i}}-\frac{W^{i}(G^{2})}{\rho_{1}^{i}})<0$, or equivalently, $\sum^{\infty}_{i=1}\frac{W^{i}(G^{1})}{\rho_{1}^{i}}<
\sum^{\infty}_{i=1}\frac{W^{i}(G^{2})}{\rho_{1}^{i}}$. 
It follows that  $f(G_{1},\rho_{1})>f(G_{2},\rho_{1})$. Note that $f(G_{1},\rho_{1})=1$.
 So, $f(G_{2},\rho_{1})<1$. Recall that $f(G_{2},x)$ is strictly increasing in its domain, and $f(G_{2},\rho_{2})=1$. We must have $\rho_{1}<\rho_{2}$. This completes the proof. \hfill$\Box$

\section{Proof of theorem \ref{cycle-planar}}

For a graph $G$, let $e(G)=|E(G)|$. If $G$ is a planar graph of order $n\geq2$, it is well-known that $e(G)\leq3n-6$. Furthermore,  $e(G)\leq2n-4$ if $G$ is also bipartite. For $u\in V(G)$, let $d_{G}(u)$ denote the degree of $u$ in $G$. For $S\subseteq V(G)$, let $d_{S}(u)$ denote the number of neighbors of $u$ in $S$. Let $K_{s,t}$ be the complete bipartite graph of color orders $s$ and $t$. Recall that $\mathcal{G}_{n,\ell}$ is the set of all $C_{\ell}$-free planar graphs of order $n$. The following Lemma \ref{spanning} is a similar result as Lemma 9 of \cite{XLF}. We need some changes, since we require this result for any $5\leq \ell\leq n$. 

(We point out that the bound $n\geq1.8\times 10^{17}$ can be decreased using the current method. However, we will not do this, as it is not the main purpose of this paper.)

\begin{lem}\label{spanning}
For integers $n\geq1.8\times 10^{17}$ and $5\leq \ell\leq n$,  any graph
$G\in{\rm SPEX}(\mathcal{G}_{n,\ell})$ contains a spanning subgraph $K_{2}\vee \overline{K_{n-2}}$.
\end{lem}

\f{\bf Proof:} Let $G\in{\rm SPEX}(\mathcal{G}_{n,\ell})$. Recall that $G$ is a planar graph. We can assume that $G$ is already embedded in the plane.
Since the resulting graph obtained by adding one suitable  edge between two components of $G$, is still a planar graph without $C_{\ell}$ but has larger spectral radius by Lemma \ref{subgraph}, we see that $G$ is connected as $G\in{\rm SPEX}(\mathcal{G}_{n,\ell})$.
Let $\rho=\rho(G)$. (Let $\mathbf{x}^{T}$ denote the transpose of a vector $\mathbf{x}$.)
By the Perron-Frobenius theorem,
there exists a positive eigenvector $\mathbf{x}=(x_1,x_2,\dots,x_n)^{\mathrm{T}}$ corresponding to $\rho$ with $\max_{u\in V(G)}x_u=1$. Let $u'\in V(G)$ with $x_{u'}=1$.
Clearly, $K_{2,n-2}$ is planar and $C_{\ell}$-free as $\ell\geq5$, which implies $K_{2,n-2}\in \mathcal{G}_{n,\ell}$. Then
$$
\rho\geq\rho(K_{2,n-2})=\sqrt{2n-4}.$$

Let $M=\{u\in V(G)~|~x_{u}\geq \frac{1}{10^4}\}$. Similar to Claims 10 and 11 of \cite{XLF}, we can show the following claim.

\medskip

\f{\bf Claim 1.} $|M|\leq\frac{n}{10^4}$, and $d_G(u)\geq (x_u-\frac{8}{10^4})n$ for any $u\in M$.

\medskip

Choose a $u''\neq u'$, such that $x_{u''}=\max_{u\in V(G)- \{u'\}}x_u$. Similar to Claim 12 of \cite{XLF}, we can show the following claim.

\medskip

\f{\bf Claim 2.} $x_{u''}\geq \frac{997}{1000}$.

\medskip

Recall that $x_{u'}=1$ and $x_{u''}\geq \frac{997}{1000}$ by Claim 2.
By Claim 1, we have
$d_{G}(u')\geq \frac{999n}{1000}$ and $d_{G}(u'')\geq \frac{996n}{1000}$.
Let $R=N_{G}(u')\cap N_{G}(u'')$ and $S=V(G)-(\{u', u''\}\cup R)$. Then $|R|\geq d_{G}(u')+d_{G}(u'')-n\geq \frac{995n}{1000}$, and $|S|\leq (n-d_{G}(u'))+(n-d_{G}(u''))\leq \frac{5n}{1000}$. Thus, 
$$\frac{|R|}{|S|}\geq\frac{995}{5}=199.$$
 Similar to Claim 13 of \cite{XLF}, we can show the following claim.

\medskip

\f{\bf Claim 3.} $x_u\leq\frac{3}{100}$ for any $u\in V(G)-\{u',u''\}$.

\medskip

The following claim is also needed in our proof.

\medskip

\f{\bf Claim 4.} $x_u\geq\frac{1}{\rho}$ for any $u\in V(G)-\{u',u''\}$.

\medskip

\f{\bf Proof of Claim 4.} Suppose $x_{u}<\frac{1}{\rho}$ for some $u\in V(G)-\{u',u''\}$. Then $1>\rho x_{u}=\sum_{v\in N_{G}(u)}x_{v}$. Let $G_{1}$ be the planar graph obtained from $G$ by deleting all edges incident with $u$ and adding an edge between $u$ and $u'$. Clearly, $G_{1}$ is also $C_{\ell}$-free. Thus, $G_{1}\in\mathcal{G}_{n,\ell}$. But
\begin{align*}
\rho(G_{1})-\rho(G)&\geq\frac{\mathbf{x}^\mathrm{T}(A(G_{1})-A(G))\mathbf{x}}{\mathbf{x}^\mathrm{T}\mathbf{x}}\\
&=\frac{2}{\mathbf{x}^\mathrm{T}\mathbf{x}}(x_{u'}-\sum_{v\in N_{G}(u)}x_{v})>0,
\end{align*}
implying that $\rho(G_{1})>\rho(G)$. But this contradicts that $G\in{\rm SPEX}(\mathcal{G}_{n,\ell})$.
\hfill$\Box$ 

\medskip

\begin{figure}
  \centering
  \includegraphics[width=11cm]{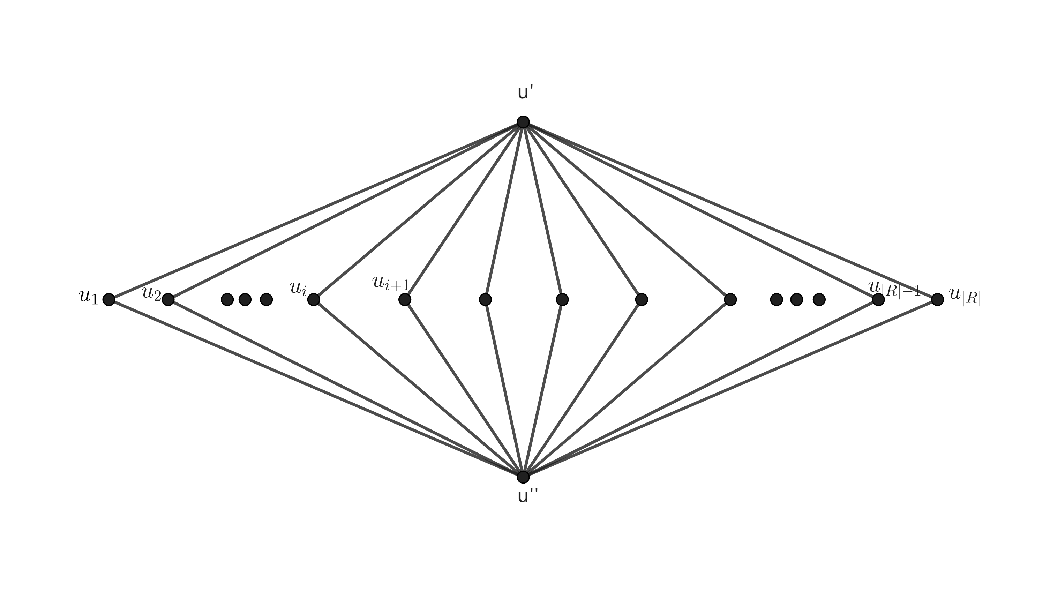}
  \caption{The local structure of the planar graph $G$}\label{figure2}
\end{figure}

In the planar embedding of $G$, we can assume that $G[\{u',u''\}\cup R]$ is embedded as in Figure \ref{figure2}, where  $u_1,u_2,\ldots,u_{|R|}\in R$ are the vertices in clockwise order around $u''$.

\medskip

\f{\bf Claim 5.} $S$ is empty.

\medskip

\f{\bf Proof of Claim 5.} Suppose to the contrary that $S$ is non-empty. Let $|S|=s\geq 1$. Note that $G$ is $K_{3,3}$-free, since $G$ is planar. Thus, for each $u\in S$, $u$ is adjacent to at most one of $u'$ and $u''$, and is adjacent to at most 2 vertices in $R$. Moreover, the maximum degree of $G[R]$ is at most 2.
Since $G[S]$ is also planar,
there exists a vertex $v_1\in S$ with $d_S(v_1)\leq 5$.
Let $S_0=S$ and $S_1=S_0- \{v_1\}$.
Repeat this process, we obtain a sequence of sets $S_0, S_1,\ldots,S_{s-1}$ such that $d_{S_{i-1}}(v_i)\leq 5$ and $S_{i}=S_{i-1}- \{v_i\}$ for each $1\leq i\leq s-1$. Let $v_{s}$ be the unique vertex in $S_{s-1}$.
Using Claim 3, for any $1\leq i\leq s$, we have
$$\sum\limits_{\substack{w\sim v_i \\ w\in \{u',u''\}\cup R\cup S_{i-1}}}x_w
\leq 1+\sum\limits_{\substack{w\sim v_i \\ w\in R }}x_w+\sum\limits_{\substack{w\sim v_i \\ w\in S_{i-1}}}x_w
\leq \frac{121}{100}<x_{u'}+x_{u''}-\frac{7}{10}.
$$
Then, noting $x_{v_{i}}\geq\frac{1}{\rho}$ by Claim 4,
$$\sum\limits_{i=1}^{s} x_{v_i}\left((x_{u'}+x_{u''})-\sum\limits_{\substack{w\sim v_i \\ w\in \{u',u''\}\cup R\cup S_{i-1}}}x_w\right)\geq\frac{7}{10}\sum\limits_{i=1}^{s} x_{v_i}\geq\frac{7}{10}\frac{s}{\rho}\geq\frac{7}{10\rho}.$$

Note that in Figure \ref{figure2}, there are $|R|-1$ closed regions: $u'u_{i}u''u_{i+1}u'$ for all $1\leq i\leq |R|-1$. Recall that $\frac{|R|}{s}\geq199$. Thus, there are at least $4$ consecutive regions containing no vertices in $S$. Without loss of generality, assume that such regions are $u'u_{i}u''u_{i+1}u'$ for $1\leq i\leq 4$. It follows that $u_{2},u_{3},u_{4}$ are not adjacent to the vertices in $S$, since $G$ is planar. Since $G$ is $K_{3,3}$-free, $u_{i}$ has at most two neighbors in $R$ for $2\leq i\leq4$.
Then,
$$\rho x_{u_i}=\sum_{u\in N_{G}(u_i)}x_u=x_{u'}+x_{u''}+\sum_{u\in N_{R}(u_i)}x_u\leq d_{G}(u_i)\leq 4$$
for each $2\leq i\leq4$.
Consequently, $x_{u_i}\leq \frac{4}{\rho}$, and hence 
$$
 x_{u_{2}}x_{u_{3}}+x_{u_{3}}x_{u_{4}}\leq \frac{32}{\rho^2}.
$$

Let $G_{2}$ be the graph obtained from $G$ by first deleting the edges $u_{2}u_{3}$, $u_{3}u_{4}$
and all the edges incident to at least one vertex in $S$, and then adding the edges $v_iu'$ and $v_iu''$  crossing the closed region $u'u_{2}u''u_{3}u'$ for each $v_i\in S$.
Clearly, $G_{2}$ is planar. Since $u_{3}$ and the vertices in $S$ are only adjacent to $u'$ and $u''$ in $G_{2}$, each $C_{k}$ with $k\geq5$ of $G_{2}$ contains at most one vertex of $S\cup \left\{u_{3}\right\}$. It follows that $G_{2}$ is also $C_{\ell}$-free, and hence $G_{2}\in \mathcal{G}_{n,\ell}$.
Recall that $x_{u_{2}}x_{u_{3}}+x_{u_{3}}x_{u_{4}}\leq \frac{32}{\rho^2}$, and
$$\sum\limits_{i=1}^{s} x_{v_i}\left((x_{u'}+x_{u''})-\sum\limits_{\substack{w\sim v_i \\ w\in \{u',u''\}\cup R\cup S_{i-1}}}x_w\right)\geq\frac{7}{10\rho}.$$
Then
\begin{equation}
\begin{aligned}
&\rho(G_{2})-\rho(G)\\
&\geq\frac{\mathbf{x}^\mathrm{T}(A(G_{2})-A(G))\mathbf{x}}{\mathbf{x}^\mathrm{T}\mathbf{x}}\nonumber\\
&\geq\frac{2}{\mathbf{x}^\mathrm{T}\mathbf{x}}
\left(\sum\limits_{i=1}^{s} x_{v_i}\left((x_{u'}+x_{u''})-\sum\limits_{\substack{w\sim v_i \\ w\in \{u',u''\}\cup R\cup S_{i-1}}}x_w\right)-x_{u_{2}}x_{u_{3}}-x_{u_{3}}x_{u_{4}}\right)\\
&\geq\frac{2}{\mathbf{x}^\mathrm{T}\mathbf{x}}(\frac{7}{10\rho}-\frac{32}{\rho^{2}})\\
&>0~~~(using~\rho\geq\sqrt{2n-4}).
\end{aligned}\notag
\end{equation}
Hence, $\rho(G_{2})>\rho(G)$,
contradicting that $G\in {\rm SPEX}(\mathcal{G}_{n,\ell})$.
Therefore, $S$ is empty. \hfill$\Box$

\medskip

\f{\bf Claim 6.} $u'u''\in E(G)$.

\medskip

\f{\bf Proof of Claim 6.} Suppose to the contrary that $u'u''\notin E(G)$. By Claim 5, we see $V(G)=\left\{u',u''\right\}\cup R$. Recall that $G[R]$ has maximum degree at most 2, as $G$ is $K_{3,3}$-free. Note that $G$ is $C_{\ell}$-free. Thus, 
$G[R]$ must be a linear forest with at least 3 components (otherwise $G$ has a Hamilton cycle).  
Then there exists some integer $1\leq i_{0}\leq n-2$ such that $u_{i_0}u_{i_0+1}\notin E(G)$.
This implies that $u'u_{i_0}u''u_{i_0+1}u'$ is a face in the embedding of $G$.

Let $G_{3}$ be the graph obtained from $G$ by adding the edge $u'u''$ crossing the face $u'u_{i_0}u''u_{i_0+1}u'$. Clearly, $G_{3}$ is planar and $C_{\ell}$-free. But $\rho(G_{3})>\rho(G)$ by Lemma \ref{subgraph}, contradicting that $G\in {\rm SPEX}(\mathcal{G}_{n,\ell})$. Therefore, $u'u''\in E(G)$. \hfill$\Box$

\medskip

From the above discussion,  we can see that $G$ contains $K_2 \vee \overline{K_{n-2}}$ as a spanning subgraph. This completes the proof. \hfill$\Box$

\medskip

The following observation is easy to see.

\medskip

\f{\bf Observation} (\cite{XLF})
Let $H= \cup_{i=1}^t P_{n_i}$, where $t\geq 2$ and $n_1\geq n_2\geq \cdots \geq n_t\geq1$. Then $K_{2}\vee H$ is $C_{\ell}$-free if and only if $n_1+n_2\leq \ell-3$.

\medskip

Now we are ready to prove Theorem \ref{cycle-planar}.

\medskip

\f{\bf Proof of Theorem \ref{cycle-planar}.} Let $G$ be a graph in ${\rm SPEX}(\mathcal{G}_{n,\ell})$. We will prove that $G$ is unique as stated in the theorem.
By Lemma \ref{spanning}, $G$ contains $K_{2}\vee \overline{K_{n-2}}$ as a spanning subgraph. Let $u',u''$ be the two dominating vertices of $G$. Let $R=V(G)-\left\{u',u''\right\}$. Since $G$ is $K_{3,3}$-free, we see that $G[R]$ has maximum degree at most 2. Since $G$ is $C_{\ell}$-free, $G[R]$ must be a linear forest with $t\geq3$ components. That is, $G[R]= \cup_{i=1}^t P_{n_i}$, where $t\geq 3,n_1\geq n_2\geq \cdots \geq n_{t}\geq1$ and $n_{1}+n_{2}\leq \ell-3$. Clearly,  $G=K_{2}\vee(\cup_{1\leq i\leq t}P_{n_{i}})$ and 
$$\sum_{1\leq i\leq t}n_{i}=n-2.$$

Now we show that $n_{1}+n_{2}=\ell-3$. Otherwise, $n_{1}+n_{2}\leq \ell-4$. Let $$G^{*}=K_{2}\vee(P_{n_{1}+1}\cup P_{n_{2}}\cup P_{n_{3}-1}\cup(\cup_{4\leq i\leq t}P_{n_{i}})).$$ Clearly, $G^{*}$ is planar and $C_{\ell}$-free (by the above observation). But $\rho(G^{*})>\rho(G)$
by Lemma \ref{1-path}, a contradiction. Hence, 
$$n_{1}+n_{2}=\ell-3.$$
Now we prove the theorem by two cases depending on the values of $\ell$.

\medskip

$(i)$ $5\leq \ell<\frac{2n+5}{3}$. Let $a\geq0,b$ be the integers, such that $n-\ell+1=a\lfloor\frac{\ell-3}{2}\rfloor+b$ with $0\leq b<\lfloor\frac{\ell-3}{2}\rfloor$.

We first show that $t\geq4$. Otherwise, $t\leq3$. Then 
$$\sum_{1\leq i\leq t}n_{i}\leq n_{1}+n_{2}+\frac{1}{2}(n_{1}+n_{2})\leq\frac{3}{2}(\ell-3)<\frac{3}{2}\frac{2n-4}{3}=n-2,$$
 a contradiction. Hence $t\geq4$. 

Now we show that $n_{1}-n_{2}\leq1$. Otherwise,  $n_{1}\geq n_{2}+2$. Let $$G^{*1}=K_{2}\vee(P_{n_{1}-1}\cup P_{n_{2}+1}\cup P_{n_{3}+1}\cup P_{n_{4}-1}\cup(\cup_{5\leq i\leq t}P_{n_{i}})).$$
 Clearly, $G^{*1}$ is planar and $C_{\ell}$-free (by the above observation). But $\rho(G^{*1})>\rho(G)$
by Lemma \ref{2-path}, a contradiction. Hence, $n_{1}-n_{2}\leq1$. It follows that $n_{1}=\lceil\frac{\ell-3}{2}\rceil$ and $n_{2}=\lfloor\frac{\ell-3}{2}\rfloor$.

Now we show that $n_{3}=n_{4}=\cdots =n_{t-1}=\lfloor\frac{\ell-3}{2}\rfloor$. Otherwise, $n_{t-1}<n_{2}=\lfloor\frac{\ell-3}{2}\rfloor$.  Let 
$$G^{*2}=K_{2}\vee((\cup_{1\leq i\leq t-2}P_{n_{i}})\cup P_{n_{t-1}+1}\cup P_{n_{t}-1}).$$
 Clearly, $G^{*2}$ is planar and $C_{\ell}$-free (by the above observation). But $\rho(G^{*2})>\rho(G)$
by Lemma \ref{1-path}, a contradiction. Hence, $n_{3}=n_{4}=\cdots n_{t-1}=\lfloor\frac{\ell-3}{2}\rfloor$. It follows that $t=a+3$ and $n_{t}=b$, since $n-\ell+1=a\lfloor\frac{\ell-3}{2}\rfloor+b$. Consequently, $G=K_{2}\vee(P_{\lceil\frac{\ell-3}{2}\rceil}\cup (a+1)P_{\lfloor\frac{\ell-3}{2}\rfloor}\cup P_{b})$, as desired.

\medskip

$(ii)$ $\frac{2n+5}{3}\leq\ell\leq n$. 
In this case, we have 
$$\ell-2+\lfloor\frac{\ell-3}{2}\rfloor\geq\frac{2n+5}{3}-2+\lfloor\frac{n-2}{3}\rfloor\geq
\frac{2n+5}{3}-2+\frac{n-4}{3}>n-2.$$

We first show that $t=3$.
Otherwise, $t\geq4$. If $n_{1}-n_{2}\geq2$, we will obtain a contradiction using a similar discussion as that in $(i)$. Thus, we must have $n_{1}-n_{2}\leq1$.  It follows that $n_{1}=\lceil\frac{\ell-3}{2}\rceil$ and $n_{2}=\lfloor\frac{\ell-3}{2}\rfloor$,
since $n_{1}+n_{2}=\ell-3$. We claim  that $n_{3}=\lfloor\frac{\ell-3}{2}\rfloor$. Otherwise, $n_{t-1}\leq n_{3}<n_{2}=\lfloor\frac{\ell-3}{2}\rfloor$.  Let 
$$G^{*3}=K_{2}\vee((\cup_{1\leq i\leq 2}P_{n_{i}})\cup P_{n_{t-1}+1}\cup P_{n_{t}-1}).$$
 Clearly, $G^{*3}$ is planar and $C_{\ell}$-free (by the above observation). But $\rho(G^{*3})>\rho(G)$
by Lemma \ref{1-path}, a contradiction. Hence, $n_{3}=\lfloor\frac{\ell-3}{2}\rfloor$. But then $$\sum_{1\leq i\leq t}n_{i}\geq n_{1}+n_{2}+n_{3}+1=\ell-2+\lfloor\frac{\ell-3}{2}\rfloor>n-2,$$
 a contradiction. Thus, we must have $t=3$. Then $G=K_{2}\vee(P_{n_{1}}\cup P_{n_{2}}\cup P_{n_{3}})$.

Now we show that $n_{2}=n_{3}$. Otherwise, $n_{2}\geq n_{3}+1$.  Let $$G^{*4}=K_{2}\vee(P_{n_{1}+1}\cup P_{n_{2}-1}\cup P_{n_{3}}).$$
 Clearly, $G^{*4}$ is planar and $C_{\ell}$-free (by the above observation). But $\rho(G^{*4})>\rho(G)$
by Lemma \ref{1-path}, a contradiction. Hence, $n_{2}=n_{3}$. It follows that $n_{1}=2\ell-n-4$ and $n_{2}=n_{3}=n-\ell+1$, since $n_{1}+n_{2}+n_{3}=n-2$ and $n_{1}+n_{2}=\ell-3$. Consequently, 
$G=K_{2}\vee(P_{2\ell-n-4}\cup2P_{n-\ell+1})$, as desired.
This completes the proof. \hfill$\Box$

\medskip

\medskip

\f{\bf Declaration of competing interest}

\medskip

There is no conflict of interest.

\medskip

\f{\bf Data availability statement}

\medskip

No data was used for the research described in the article.


\medskip

\begin{thebibliography}{99}









\bibitem{BR}
 B.N. Boots and G.F. Royle, A conjecture on the maximum value of the principal eigenvalue of a planar graph, Geogr. Anal. 23 (3) (1991) 276-282.

\bibitem{CRS}
D. Cvetkovi\'{c}, P. Rowlinson and S. Simi\'{c}, An Introduction to the Theory of Graph Spectra, Cambridge University Press, Cambridge, 2010.

\bibitem{CV}
 D. Cao and A. Vince, The spectral radius of a planar graph, Linear Algebra Appl. 187 (1993)
251-257.




\bibitem{FLS}
L. Fang, H. Lin and Y. Shi, Extremal spectral results of planar graphs without vertex-disjoint cycles. J. Graph Theory 106 (2024) 496-524.






\bibitem{N}
V. Nikiforov, Bounds on graph eigenvalues II, Linear Algebra Appl. 427 (2007) 183-189.


\bibitem{TT}
M. Tait and J. Tobin, Three conjectures in extremal spectral graph theory, J. Combin. Theory, Ser. B 126
(2017) 137-161.


\bibitem{XLF}
P. Xu, H. Lin and L. Fang, Long cycles and spectral radii in planar graphs, Electron. J. Comb. 32 (2025) \#P2.26.


\bibitem{ZW}
M. Zhai and B. Wang, Proof of a conjecture on the spectral radius of $C_{4}$-free graphs,
Linear Algebra Appl. 437 (2012) 1641-1647.


\bibitem{Z}
W. Zhang, Walks, infinite series and spectral radius of graphs, arXiv:2406.07821v3.

\end{thebibliography}
\end{document}